\documentstyle[12 pt]{article}
\def\endpf{\hbox{\vrule height1.5ex width.5em}}

\def\d{\displaystyle}
\begin{document}
\title{Stabilization with  finite dimensional controllers for a periodic parabolic system
under perturbations in the system conductivity}
\author{  Ling Lei \footnote{
Supported  by a National Science Foundation of China Research Grant (NSFC-10801108)}\\  Department of Mathematics and Statistics, Wuhan University,\\
 Wuhan, 430072, P.R.China}
\date{}
\maketitle
{\bf Abstract.}    This work studies the stabilization for a
periodic parabolic system under perturbations in the system
conductivity.   A  perturbed system does not have any periodic
solution in general. However, we will prove that the perturbed
system can always be pulled back to a periodic system after imposing
a control from a fixed finite dimensional subspace.
The paper continues the author's previous work in
\cite{kn:[1]}.

\bigskip

\vskip 1cm

\section{Introduction}
\hspace*{0.5 cm}Let $\Omega\subset {\bf R}^N$ be a bounded domain
with a $C^2$-smooth boundary $\partial\Omega$ and let
$\omega\subset\Omega$ be a subdomain. Write $Q=\Omega \times (0,T)$
with $T>0$ and  write $\Sigma=\partial\Omega\times (0,T) $. Consider
the following parabolic equation:
$$
 \left\{\begin{array}{ll}
\displaystyle{\frac{\partial y}{\partial t}(x,t)}+L_0
y(x,t)+e(x,t)y(x,t)=f(x,t), \;&
\mbox{in }\;Q=\Omega \times (0,T),\\
 y(x,t)=0, & \mbox{on }\;\Sigma=\partial \Omega \times (0,T),\\
\end{array}\right.
\eqno{(1.1)}
$$
where
$$
L_0y(x,t) = - \sum ^N _{i,j=1} \frac{\partial }{\partial x_j}
(a^{ij}(x)\frac{\partial }{\partial x_i} y(x,t)) +c(x)y(x,t)
$$
is considered as the system operator. Here and in all that  follows,
we make the following regularity assumptions for the coefficients of
$L_0$:

\noindent (I):
$$\begin{array}{ll}
a^{ij}(x) \in Lip(\overline{\Omega}),\;a^{ij}(x)=a^{ji}(x),\;
\mbox{and}\; \lambda ^*|\xi|^2 \leq \d{\sum_{i,j=1}^{N}} a^{ij}(x)
\xi _i \xi_j & \leq \displaystyle{\frac{1}{\lambda^*}} |\xi|^2 ,\;
\mbox{for}\; \xi
\in {\bf R}^N \\
\end{array}
\eqno{(1.2)}
$$
 with $\lambda^*$  a  certain positive constant;

\noindent (II):  $$
\begin{array}{ll}
 c(x) \in
L^\infty(\Omega),\;
 e(x,t) \in L^\infty (0,T;L^q(\Omega))\ \hbox{with }\ q
>\max\{N,2\},\  \hbox{and}\  f(x,t)\in L^2(Q).
\end{array}
\eqno{(1.3)}
$$
In such a system, we regard $e(x,t)$ as a perturbation in the system
conductivity. Suppose in the ideal case, namely, in the case when
the perturbation $e(x,t)\equiv 0$, (1.1) has a periodic solution
$y_0(x,t)$:
$$
 \left\{\begin{array}{ll}
\displaystyle{\frac{\partial y_0}{\partial t}(x,t)}+L_0
y_0(x,t)=f(x,t), \;\;\;\;\;\;\;\;\;\;\;\;\;\;\;\;\;\;\;\;\;&
\mbox{in }\;\;Q,\\
 y_0(x,t)=0, & \mbox{on }\;\; \Sigma,\\
y_0(x,0)=y_0(x,T),&\mbox{in }\Omega.
\end{array}\right.
\eqno{(1.4)}
$$
Then the presence of the error term $e(x,t)$ may well destroy the
periodicity of the system. Indeed, (1.1) may no longer have any
periodic solution. (See Section 3.) The problem that we are
interested in in this paper is to understand if there is a finite
(constructible) dimensional subspace ${\bf U}\subset L^2(Q)$, such
that, after imposing a control $u_e\in \bf{U}$, we can restore the
periodic solution $y_e$. Moreover, we would like to know if $y_e$ is
close to $y_0$ and if the energy of $u_e$ is small, when $e(x,t)$ is
small. Our main purpose of this paper is to show that we can indeed
achieve this goal in the small perturbation case, even if the
control is only imposed over a subregion $\omega$ of $\Omega$. The
basic tool for this study is the existence and energy estimate for
the approximate periodic solutions obtained in the author's previous
paper \cite{kn:[1]}. \medskip

To state our results, we first recall the definition of approximate
periodic solutions with respect to the elliptic operator $L_0$.

Notice that $L_0$ is a symmetric operator. Consider the eigenvalue
problem of $L_0$:
$$
\left\{
\begin{array}{ll}
L_0 X(x)=\lambda X(x),\\
X(x)|_{\partial \Omega} =0.
\end{array}\right.
\eqno{(1.5)}
$$
Making use of the regularity assumptions of the coefficients of
$L_0$, we know (see, \cite{kn:[2]} \cite{kn:[3]}, for example) that
(1.5) has a complete set of eigenvalues
$\{\lambda_j\}_{j=1}^{\infty}$ with the associated eigenvectors
$\{X_j(x)\}_{j=1}^{\infty}$ such that $$L_0 X_j(x)=\lambda_j X_j
(x),$$ $$-\infty<\lambda_1\leq\lambda_2
\leq\cdots\leq\lambda_j\leq\cdots<\infty,\ \lim_{j\rightarrow
\infty}\lambda_j=\infty,\;X_j(x)\in H^1_0(\Omega)\cap
C(\overline{\Omega}).$$ Choose $\{X_j(x)\}_{j=1}^{\infty}$ such that
it forms an orthonormal basis of $L^2(\Omega)$. Therefore, for any
$y(x,t)\in L^2(Q)$, we have $
y(x,t)=\d\sum^{\infty}_{j=1}y_j(t)X_j(x)$, where
$$y_j(t)=\langle y(x,t),X_j(x)\rangle= \d\int_\Omega
y(x,t)X_j(x)dx\in L^2(0,T).$$\medskip

{\bf Definition 1.1.} {\it We call $y(x,t)$ is a K-approximate
periodic solution of (1.1) with respect
to $L_0$ if \\
(a):  $y \in C([0,T];L^2(\Omega))\cap L^2(0,T;H^1_0 (\Omega ))$ is
a weak solution of (1.1);\\
(b): $ y\in {\bf S_{K}} $, where ${\bf S_{K}}$ is the space of the
following functions: $${\bf S_{K}} = \{y(x,t)\in
L^2(Q);\;y_j(0)=y_j(T),\;\mbox{for}\;j \ge K+1,\;
y_j(t)=\displaystyle{\int_\Omega } y(x,t)X_j (x)dx \}.$$}

\medskip
When $K= 0$, we will always  regard $\sum ^{0} _{j=1} = 0$. Hence, a
0-approximate periodic solution of (1.1) is a regular periodic
solution. In what follows, we write $\langle y(\cdot,t),y(\cdot
,t)\rangle=\displaystyle{\int_\Omega} y^2(x,t)dx=\|y(\cdot,t)\|^2$,
and we denote $y_t$ for the derivative of $y(x,t)$ with respect to
$t$.

\medskip
Our first result of this paper can be stated as follows:

\bigskip
{\bf Theorem 1.1} {\it  Consider the system (1.1), where $e(x,t)$ is
regarded as a perturbation in the system conductivity. Suppose that
(1.1) has a periodic solution $y_0(x,t)$ at the ideal case with
$e(x,t)\equiv 0$. Assume that
$\|e(x,t)\|_{L^\infty(0,T;L^q(\Omega))}=ess\d{\sup_{t\in(0,T)}}\|e(\cdot,t)\|_{L^q(\Omega)}<\varepsilon$,
where $\varepsilon<1$ is a small constant which depends only on
$L_0,\Omega, N, q, T$ with $q>\max\{N,2\}$. Then there are a
non-negative integer $K_0$, depending only on $L_0,\Omega, N, q, T$
(but not $f$), and  a unique outside force of the form
$$u_e(x):=\sum_{j=1}^{K_0}u_jX_j(x)\in {\bf U}=span_{{\mathbf
R}}\{X_1(x),X_2(x),\cdots,X_{K_0}(x)\},$$ where $u_j\in {\mathbf
R}$, such that the following has a unique periodic solution $y$
satisfying:
$$
 \left\{\begin{array}{ll}
\displaystyle{\frac{\partial y(x,t)}{\partial
t}}+L_0y(x,t)+e(x,t)y(x,t)=f(x,t)+u_e(x),
\;\;\;\;\;\;\;\;\;\;\;\;\;\;\;\;\;\;\;\;\;&
\mbox{in }\;\;Q,\\
 y(x,t)=0, & \mbox{on }\;\; \Sigma,\\
 \langle y(x,0),X_j(x)\rangle=\langle y_0(x,0),X_j(x)\rangle, & \mbox{for }\;j\leq K_0,\\
y(x,0)=y(x,T),&\mbox{in }\;\;\Omega.
\end{array}\right.
\eqno{(1.6)}
$$
Moreover,  we have the following energy estimate:
$$\begin{array}{ll}
&\d\sup_{t\in[0,T]}\|(y-y_0)(\cdot,t)\|^2+\displaystyle{\int^T_0}\|\nabla(y-y_0)(\cdot,t)\|^2dt\\
&\leq C(system,K_0)\|e(x,t)\|^2_{L^\infty(0,T;L^q(\Omega))}(1
+|\vec{a}|^2+\d{\int_Q} f^2dxdt),\end{array}\eqno{(1.7)}
$$
and
$$
\|u_e\|^2_{L^2(\Omega)}\leq
C(system,K_0)\|e(x,t)\|^2_{L^\infty(0,T;L^q(\Omega))}(1+|\vec{a}|^2+\d{\int_Q}
f^2dxdt), \eqno{(1.8)}
$$
where $\vec{a}=(a_1,a_2,\cdots,a_{K_0})=(\langle
y_0(x,0),X_1(x)\rangle,\langle y_0(x,0),X_2(x)\rangle,\cdots,\langle
y_0(x,0),X_{K_0}(x)\rangle)$.
Here and in what follows, $C(system,K_0)$
denotes a constant depending only on $L_0,\Omega, N,q,T$, which may
be different in different contexts.}\bigskip


\bigskip
In Section 3 of this paper, we will construct an example, showing
that without outside controls, (1.1)  has no periodic solutions in
general. This is one of the main features in our  Theorem 1.1: The
control can always be taken from a certain fixed constructible {\it
finite dimensional subspace} to regain the periodicity, while
 the perturbation space for $e(x,t)$, which destroys the periodicity, is
{\it of infinite dimension}. We also notice  that our system
operator $L_0$ is not assumed to be positive.

\medskip

The second part of this work is to consider the same problem as
studied in the first part, but with the control only imposed over a
subregion $\omega\subset\Omega$ and time interval $E\subset [0,T]$,
$m(E)>0$. We will similarly obtain  the following:\bigskip

{\bf Theorem 1.2.} {\it Suppose that the system (1.1) has a periodic
solution $y_0(x,t)$ at the ideal case with $e(x,t)\equiv 0$. Then
there are a positive integer $K_0$, a small constant
$\varepsilon>0$, depending only on $L_0,\Omega, N,q,T$
$(q>\max\{N,2\})$, such that, when
$$\|e(x,t)\|_{L^\infty(0,T;L^q(\Omega))}=ess\d{\sup_{t\in(0,T)}}\|e(x,t)\|_{L^q(\Omega)}<\varepsilon,$$
the following has a unique periodic solution:
$$
 \left\{\begin{array}{ll}
\displaystyle{\frac{\partial y(x,t)}{\partial
t}}+L_0y(x,t)+e(x,t)y(x,t)=f(x,t)+\displaystyle{\sum_{j=1}^{K_0}}\chi_\omega(x)\chi_E(t)u_jX_j(x),
\;\;\;\;\;\;\;\;&
\mbox{in }\;\;Q,\\
 y(x,t)=0, & \mbox{on }\;\; \Sigma,\\
 \langle y(x,0),X_j(x)\rangle=a_j, & \mbox{for }\;j\leq K_0,\\
y\in {\bf S_{K_0}},
\end{array}\right.
\eqno{(1.9)}
$$
where $(a_1,a_2,\cdots,a_{K_0})=(\langle
y_0(x,0),X_1(x)\rangle,\langle y_0(x,0),X_2(x)\rangle,\cdots,\langle
y_0(x,0),X_{K_0}(x)\rangle)=\vec{a}$,
$(u_1,u_2,\cdots,u_{K_0})=\vec{u}\in {\mathbf R}^{K_0}$. Moreover,
$$
|\vec{u}|^2\leq
C(system,K_0,\omega)\d\frac{\|e(x,t)\|^2_{L^\infty(0,T;L^q(\Omega))}}{(m(E))^2}(1+|\vec{a}|^2+\d{\int_Q}
f^2dxdt),\eqno{(1.10)}
$$
and
$$\begin{array}{ll}
&\d\sup_{t\in[0,T]}\|(y-y_0)(\cdot,t)\|^2+\displaystyle{\int^T_0}\|\nabla(y-y_0)(\cdot,t)\|^2dt\\
&\leq
C(system,K_0,\omega)\d\frac{\|e(x,t)\|^2_{L^\infty(0,T;L^q(\Omega))}}{(m(E))^2}(1+|\vec{a}|^2+\d{\int_Q}
f^2dxdt).\end{array}\eqno{(1.11)}
$$
Here, $$\chi_\omega(x),\;\chi_E(t)$$ are the characteristic
functions for  $\omega$ and $E$, respectively; and $C(system,K_0,
\omega)$ is a constant depending only on $\omega,L_0,\Omega, N,q,T$.
}

\bigskip

Theorem 1.1 and Theorem 1.2 give  stabilization results for the
periodic solutions of a linear parabolic system under small
perturbation of the system conductivity, modifying a control from a
fixed finite dimensional subspace. We do not know if  similar
results as in Theorem 1.1 hold under the large perturbation case.

\medskip
The paper is organized as follows. In Section 2, we prove Theorem
1.1. In Section 3, we give an example to show that with a small
perturbation $e(x,t)$, (1.1) has no periodic solution in general. In
section 4, we give the proof of Theorem 1.2.
\section{Small perturbation}
 \hspace*{0.5 cm}In this Section, we give a proof of Theorem
 1.1, based on the author's previous paper \cite{kn:[1]}. For  convenience of the reader, we first recall the following result of   \cite{kn:[1]}, which will be used here. 
\medskip

{\bf Theorem 2.1}. {\it Assume (1.2) and (1.3). Let $e(x,t)\in {\mathcal{M}}(q,M)$, where, for any positive
number $M$ and $q> \frac{N}{2}$,
 $${\mathcal{M}}(q,M):=
\{e(x,t) \in L^{\infty}(0,T;L^q (\Omega)); \mbox{ess sup}_{t\in
(0,T)} \|e(x,t)\|_{L^q (\Omega)}\le M\}.$$ Then, there exists an
integer $K_0(L_0,M,\Omega, q,N,T)$ $\geq 0$, depending only on
$(L_0, M,\Omega, q,N, T)$ (but not $f(x,t)$), such that for any
$K\geq K_0(L_0,M,\Omega,q,N, T)$ and any initial value
$\vec{a}=(a_1,a_2,\cdots,a_K)\in {\bf R^K}$, we have a unique
solution to the following equation:
$$
 \left\{\begin{array}{ll}
\displaystyle{\frac{\partial y(x,t)}{\partial
t}}+L_0y(x,t)+e(x,t)y(x,t)=f(x,t),
\;\;\;\;\;\;\;\;\;\;\;\;\;\;\;\;\;\;\;\;\;&
\mbox{in }\;\;Q,\\
 y(x,t)=0, & \mbox{on }\;\; \Sigma,\\
 \langle y(x,0),X_j(x)\rangle=a_j, & \mbox{for }\;j\leq K,\\
y\in {\bf S_{K}}.
\end{array}\right.
\eqno{(2.1)}
$$
Moreover, for such a solution $y(x,t)$, we have the following energy estimate:
$$
\begin{array}{ll}
\displaystyle{\sup_{t\in[0,T]}}\|y(\cdot,t)\|^2
+\displaystyle{\int^T_0}\|\nabla y(\cdot,t)\|^2 dt \leq
C(L_0,M,\Omega,q,N,T) (|\vec{a}|^2 +\displaystyle{\int_Q}f^2dxdt).
\end{array}
\eqno (2.2) $$}
\bigskip

Now, suppose $y_0$ is a periodic solution of (1.1) with $e(x,t)=0$,
namely,
$$
 \left\{\begin{array}{ll}
\displaystyle{\frac{\partial y_0}{\partial t}(x,t)}+L_0
y_0(x,t)=f(x,t), \;\;\;\;\;\;\;\;\;\;\;\;\;\;\;\;\;\;\;\;\;&
\mbox{in }\;\;Q,\\
 y_0(x,t)=0, & \mbox{on }\;\; \Sigma,\\
y_0(x,0)=y_0(x,T),&\mbox{in }\Omega.
\end{array}\right.
\eqno{(2.3)}
$$
Let $$ a_j=(y_0)_j(0)=\langle y_0(x,0),X_j(x)\rangle,\mbox{ for }
j=1,2,\cdots.
$$
In all that follows, we assume that  $e(x,t)\in {\mathcal{M}}(q,M)$
with $M=1$. By Theorem 2.1, there exists an integer
$K_0(L_0,M,\Omega, q,N,T)$ $\geq 0$, such that for the initial value
$\vec{a}=(a_1,a_2,\cdots,a_{K_0})\in {\mathbf R}^{K_0}$, we have a
unique solution $y(x,t)$ satisfying the following equations:
$$
 \left\{\begin{array}{ll}
\displaystyle{\frac{\partial y(x,t)}{\partial
t}}+L_0y(x,t)+e(x,t)y(x,t)=f(x,t)+\displaystyle{\sum_{j=1}^{K_0}}u_jX_j(x),
\;\;\;\;\;\;\;\;&
\mbox{in }\;\;Q,\\
 y(x,t)=0, & \mbox{on }\;\; \Sigma,\\
 \langle y(x,0),X_j(x)\rangle=a_j, & \mbox{for }\;j\leq K_0,\\
y\in {\bf S_{K_0}}.
\end{array}\right.
\eqno{(2.4)}
$$Here, $\vec{u}=(u_1,u_2,\cdots,u_{K_0})\in {\mathbf R}^{K_0}$.

Subtracting (2.3) from (2.4), we get the following equation:
$$
 \left\{\begin{array}{ll}
(y-y_0)_t
+L_0(y-y_0)+e(x,t)(y-y_0)=\displaystyle{\sum_{j=1}^{K_0}}u_jX_j(x)-e(x,t)y_0,
\;\;\;&
\mbox{in }\;\;Q,\\
 (y(x,t)-y_0(x,t))=0, & \mbox{on }\;\; \Sigma,\\
 (y-y_0)_j(0)=\langle y(x,0)-y_0(x,0),X_j(x)\rangle=0, & \mbox{for }\;j\leq K_0,\\
(y-y_0)\in {\bf S_{K_0}}.
\end{array}\right.
\eqno{(2.5)}
$$
We define a map$$J:\;{\mathbf R}^{K_0}\longmapsto {\mathbf
R}^{K_0}$$by
$$J(u_1,u_2,\cdots,u_{K_0})=((y-y_0)_1(T),(y-y_0)_2(T),\cdots,(y-y_0)_{K_0}(T)).$$Write $v=y-y_0=v_0+v_u$. Here,
$v_0$ and $v_u$ are the solution of the following equations, respectively,
$$
 \left\{\begin{array}{ll}
(v_0)_t +L_0v_0+e(x,t)v_0=-e(x,t)y_0, \;\;\;&
\mbox{in }\;\;Q,\\
 v_0=0, & \mbox{on }\;\; \Sigma,\\
(v_0)_{j}(0)=0, & \mbox{for }\;j\leq K_0,\\
v_0\in {\bf S_{K_0}}.
\end{array}\right.
\eqno{(2.6)}
$$and
$$
\left\{\begin{array}{ll} (v_u)_t
+L_0v_u+e(x,t)v_u=\displaystyle{\sum_{j=1}^{K_0}}u_jX_j(x), \;\;\;&
\mbox{in }\;\;Q,\\
 v_u=0, & \mbox{on }\;\; \Sigma,\\
(v_u)_{j}(0)=0, & \mbox{for }\;j\leq K_0,\\
v_u\in {\bf S_{K_0}}.
\end{array}\right.
\eqno{(2.7)}
$$
We are led to the question to find out if there is a vector $\vec{u}=(u_1,u_2,\cdots,u_{K_0})\in {\bf R^{K_0}}$
such that
$$
J(\vec{u})=((y-y_0)_1(T),(y-y_0)_2(T),\cdots,(y-y_0)_{K_0}(T))=(0,0,\cdots,0).
$$
Indeed, if this is the case, then $y$ is a periodic solution with
the required estimate as we will see later.

 For this purpose, we write
$J_0=((v_0)_1(T),(v_0)_2(T),\cdots,(v_0)_{K_0}(T))$ and
$$J^*(\vec{u})=((v_u)_1(T),(v_u)_2(T),\cdots,(v_u)_{K_0}(T)).$$Then $$J(\vec{u})=J_0+J^*(\vec{u}).$$
Now, it is easy to see that $J^*$ is linear in
$(u_1,u_2,\cdots,u_{K_0})$. We next claim that $J^*$ is invertible
under the small perturbation case. If not, we can find a vector
$\vec{\xi}=(\xi_1,\xi_2,\cdots,\xi_{K_0})\in {\mathbf R}^{K_0}$ with
$|\vec{\xi}|=\sqrt{\xi^2_1+\xi^2_2+\cdots+\xi^2_{K_0}}=1$ such that
$J^*(\vec{\xi})=0$. Hence, we have a unique solution to the
following problem:
$$
\left\{\begin{array}{ll} w_t
+L_0w+e(x,t)w=\displaystyle{\sum_{j=1}^{K_0}}\xi_jX_j(x), \;\;\;&
\mbox{in }\;\;Q,\\
w=0, & \mbox{on }\;\; \Sigma,\\
w_{j}(0)=w_{j}(T)=0, & \mbox{for }\;j\leq K_0,\\
w\in {\bf S_{K_0}}.
\end{array}\right.
\eqno{(2.8)}
$$
First, by the energy estimate in Theorem 2.1, we have for $w(x,t)$,
$$\begin{array}{ll}
\displaystyle{\sup_{t\in[0,T]}}\|w(\cdot,t)\|^2
+\displaystyle{\int^T_0}\|\nabla w(\cdot,t)\|^2 dt &\leq
C(system)\cdot T\cdot |\vec{\xi}|^2\\ &\leq C(system,K_0).
\end{array}
\eqno (2.9)$$As mentioned before, we use $C(system, K_0)$ to denote
a constant depending only on $L_0,M,\Omega,q,N,T$, which may be
different in different contexts.
\medskip

Write $w=\displaystyle{\sum^{\infty}_{j=1}}w_j(t)X_j(x)$ as before. Then we have
$$
\displaystyle{\frac{dw_j(t)}{dt}}+\lambda_j w_j(t)+\int_\Omega
e(x,t)w(x,t)X_j(x)dx =\xi_j,\;\mbox{for
}j=1,2,\cdots,K_0.\eqno{(2.10)}
$$
Next, by the H$\ddot{o}$lder inequality (see Claim 2.2 of
\cite{kn:[1]}), we have
$$\begin{array}{ll}
\displaystyle{\int_\Omega}|e(x,t)w(x,t)X_j(x)|dx&\leq
C(\Omega,N,q)\|e(x,t)\|_{L^\infty(0,T;L^q(\Omega))}[\|w(\cdot,t)\|^2_{L^2(\Omega)}\\
&+\|X_j(x)\|^2_{L^2(\Omega)}+ \|\nabla w(\cdot,t)\|^2_{L^2(\Omega)}+\|\nabla X_j(x)\|^2_{L^2(\Omega)}]\\
&\leq C(\Omega,N,q)\|e(x,t)\|_{L^\infty(0,T;L^q(\Omega))}[1+\lambda^2_j\\
&+ \|w(\cdot,t)\|^2_{L^2(\Omega)}+
\|\nabla w(\cdot,t)\|^2_{L^2(\Omega)}].
\end{array}
$$
By (2.9), we have
$$
\displaystyle{\int^T_0\int_\Omega}|e(x,t)w(x,t)X_j(x)|dxdt\leq
C(system,K_0)\|e(x,t)\|_{L^\infty(0,T;L^q(\Omega))}.\eqno{(2.11)}
$$
Next, from (2.10), we get
$$(e^{\lambda_j t}w_j(t))'_t+\int_\Omega e(x,t)w(x,t)X_j(x)e^{\lambda_j t}dx
=e^{\lambda_j t}\xi_j,\;\mbox{for }j=1,2,\cdots,K_0.$$ Integrating the above over [0,T], we get, for
$j=1,2,\cdots,K_0$,
$$0+\displaystyle{\int^T_0\int_\Omega}e(x,t)w(x,t)X_j(x)e^{\lambda_j t}dxdt=\xi_j\d{\int^T_0}e^{\lambda_j t}dt.$$
Namely,
$$
\xi_j=\left\{\begin{array}{ll} \d{\frac{\displaystyle{\int^T_0\int_\Omega}e(x,t)w(x,t)X_j(x)e^{\lambda_j
t}dxdt}{\frac{1}{\lambda_j}(e^{\lambda_j T}-1)}},\;\;\;\;&\mbox{for }\lambda_j\neq 0,\\
\d{\frac{\displaystyle{\int^T_0\int_\Omega}e(x,t)w(x,t)X_j(x)dxdt}{T}},&\mbox{for }\lambda_j=0.
\end{array}\right.
$$
Hence, we get, for $j=1,2,\cdots,K_0$,
$$\begin{array}{ll}
|\xi_j|&\leq C(system,K_0)\displaystyle{\int^T_0\int_\Omega}|e(x,t)w(x,t)X_j(x)|dxdt\\
&\leq C(system,K_0)\|e(x,t)\|_{L^\infty(0,T;L^q(\Omega))}.
\end{array}\eqno{(2.12)}
$$
We get
$$1=|\vec{\xi}|\leq C(system,K_0)\sqrt{K_0}\|e(x,t)\|_{L^\infty(0,T;L^q(\Omega))}.$$
This gives a contradiction when
$$\|e(x,t)\|_{L^\infty(0,T;L^q(\Omega))}<
\d{\frac{1}{C(system,K_0)\sqrt{K_0}}}.$$ Therefore, we showed that
$J^*$ is invertible when
$\|e(x,t)\|_{L^\infty(0,T;L^q(\Omega))}<\epsilon $ with a certain
$\epsilon$ depending only on $L_0,\Omega, N,q, T$.

 Hence, for any given
$\vec{b}=(b_1,b_2,\cdots,b_{K_0})\in {\mathbf R}^{K_0}$, there
exists a unique $$\vec{u}=(u_1,u_2,\cdots,u_{K_0})\in {\mathbf
R}^{K_0}$$ such that
$$J^*(\vec{u})=J^*(u_1,u_2,\cdots,u_{K_0})=(b_1,b_2,\cdots,b_{K_0}).$$
Back to the equation (2.7), we have
$$
\displaystyle{\frac{d(v_u)_j(t)}{dt}}+\lambda_j (v_u)_j(t)+\int_\Omega e(x,t)v_u(x,t)X_j(x)dx =u_j,\;\mbox{for
}j=1,2,\cdots,K_0.
$$
Then
$$
\displaystyle{\frac{d[e^{\lambda_j t}(v_u)_j(t)]}{dt}}+\int_\Omega e(x,t)v_u(x,t)X_j(x)e^{\lambda_j t}dx=u_j
e^{\lambda_j t},\;\mbox{for }j=1,2,\cdots,K_0.
$$
Integrating the above over [0,T], by the definition of $J^*$, we have
$$
b_j e^{\lambda_j T}-0+\d{\int^T_0\int_\Omega} e(x,t)v_u(x,t)X_j(x)e^{\lambda_j t}dxdt=u_j \int^T_0e^{\lambda_j
t}dt,\;\mbox{for }j=1,2,\cdots,K_0.
$$
We then get
$$
u_j=\left\{\begin{array}{ll} \d{\frac{b_j e^{\lambda_j
T}+\displaystyle{\int^T_0\int_\Omega}e(x,t)v_u(x,t)X_j(x)e^{\lambda_j
t}dxdt}{\frac{1}{\lambda_j}(e^{\lambda_j T}-1)}},\;\;\;\;&\mbox{for }\lambda_j\neq 0,\\
\d{\frac{\displaystyle{\int^T_0\int_\Omega}e(x,t)v_u(x,t)X_j(x)dxdt}{T}},&\mbox{for }\lambda_j=0.
\end{array}\right.
$$
$$\begin{array}{ll}
|u_j|^2 &\leq 2e^{2\lambda_{K_0}T}|b_j|^2
+2e^{2\lambda_{K_0}T}[\displaystyle{\int^T_0\int_\Omega}e(x,t)v_u(x,t)X_j(x)dxdt]^2\\
&\leq 2e^{2\lambda_{K_0}T}|b_j|^2 +2e^{2\lambda_{K_0}T}\cdot
\hbox{sup}_{\Omega}|X_j|^2[\displaystyle{\int^T_0}\|e(\cdot,t)\|_{L^q(\Omega)}\|v_u(\cdot,t)\|_{L^{q'}(\Omega)}dt]^2
\end{array}
$$
Here $1/q+1/q'=1$. Since $\Omega$ is bounded and
$q'=\frac{q}{q-1}\le 2$, by the H\"older inequality, we have
$\|v_u\|_{L^{q'}(\Omega)}\le C(\Omega, q)\|v_u\|_{L^{2}(\Omega)}.$
Hence,
$$[\displaystyle{\int^T_0}\|e(\cdot,t)\|_{L^q(\Omega)}\|v_u(\cdot,t)\|_{L^{q'}(\Omega)}dt]^2\le
C(\Omega,T,q)\|e\|^2_{L^\infty(0,T;L^q(\Omega))}\|v_u\|^2_{L^{2}(Q)}.$$

By the energy estimate in Theorem 2.1, we have
$\|v_u\|^2_{L^{2}(Q)}\le C(system, K_0) |u|^2.$ Hence, as argument
before, when $\|e\|^2_{L^\infty(0,T;L^q(\Omega))}$ is small, we can
solve the above to obtain the following:

$$|\vec{u}|^2 \leq C(system,K_0)|\vec{b}|^2.\eqno{(2.13)}$$

Back to (2.5), we need to find $\vec{u}=(u_1,u_2,\cdots,u_{K_0})$
such that the solution in (2.5) has the property $(y-y_0)_j(T)=0$
for $j=1,2,\cdots,K_0$. As mentioned before,  $v=y-y_0$ is then a
periodic solution. Thus $y=v+y_0$ is  a periodic solution of (2.4)
after applying the control force
$\displaystyle{\sum_{j=1}^{K_0}}u_jX_j(x)$. To this aim, we need
only to find $\vec{u}$ such that
$$J(\vec{u})=0\;\mbox{or }\;J^*(\vec{u})=-J_0.$$ By the definition
of $J_0$, $J_0=-\vec{b}=(-b_1,-b_2,\cdots,-b_{K_0})$ is given by
$$
 \left\{\begin{array}{ll}
(v_0)_t +Lv_0+e(x,t)v_0=-e(x,t)y_0, \;\;\;&
\mbox{in }\;\;Q,\\
 v_0=0, & \mbox{on }\;\; \Sigma,\\
(v_0)_{j}(0)=0,\;(v_0)_j(T)=-b_j, & \mbox{for }\;j\leq K_0,\\
v_0\in {\bf S_{K_0}}.
\end{array}\right.
$$
By the energy estimate of Theorem 2.1, we have
$$\begin{array}{ll}
|\vec{b}|^2&\leq \|v_0(\cdot,T)\|^2_{L^2(\Omega)}\\
&\leq C(system,K_0)\d{\int^T_0\int_\Omega (-ey_0)^2dxdt}\\
&\leq C(system,K_0)\d{\int^T_0}\{\|e(x,t)\|^2_{L^q(\Omega))}\|y_0(\cdot,t)\|^2_{L^{\frac{2q}{q-2}}(\Omega)}\}dt\\
&\leq C(system,K_0)\|e(x,t)\|^2_{L^\infty(0,T;L^q(\Omega))}\|\nabla y_0\|^2_{L^2(Q)}\\
&\leq
C(system,K_0)\|e(x,t)\|^2_{L^\infty(0,T;L^q(\Omega))}(|\vec{a}|^2+\d{\int_Q}
f^2dxdt),
\end{array}\eqno{(2.14)}
$$
where $\vec{a}=(a_1,a_2,\cdots,a_{K_0})=(\langle y_0(x,0),X_1(x)\rangle,\langle
y_0(x,0),X_2(x)\rangle,\cdots,\langle y_0(x,0),X_{K_0}(x)\rangle)$.

Thus, by (2.13), we get
$$|\vec{u}|^2\leq C(system,K_0)\|e(x,t)\|^2_{L^\infty(0,T;L^q(\Omega))}
(1+|\vec{a}|^2+\d{\int_Q} f^2dxdt).\eqno{(2.15)}
$$
By (2.2), (2.14) and (2.15), we obtain
$$\begin{array}{ll}
&\d\sup_{t\in[0,T]}\|(y-y_0)(\cdot,t)\|^2+\displaystyle{\int^T_0}\|\nabla(y-y_0)(\cdot,t)\|^2dt\\
&\leq C(system,K_0)\|e(x,t)\|^2_{L^\infty(0,T;L^q(\Omega))}(1
+|\vec{a}|^2+\d{\int_Q} f^2dxdt),\end{array}
$$
Summarizing the above, we complete the proof of Theorem 1.1.
$\endpf$
\bigskip


\section{An example}
\hspace*{0.5 cm}In this section, we present an example, showing that
with a small perturbation $e(x,t)$, (1.1) has no periodic solution
in general. This demonstrates the importance of an outside control
to gain back the periodicity as in Theorem 1.1.\medskip

We consider the following one dimensional parabolic equation:
$$
\left\{\begin{array}{ll} y_t-y_{xx}-y-e(x)y=f(x),\;\;\;\;\;\;&0\leq x\leq \pi,\;0\leq t\leq T,\\
y(0,t)=y(\pi,t)=0,&0\leq t\leq T.
\end{array}\right.\eqno{(3.1)}
$$
Let $L_e y=-y_{xx}-y-e(x)y$ with $e(x)\in C^0[0,\pi]$. Suppose $0$
is an eigenvalue of $L_e$ with eigenvectors $\{X_j(x)\}^m_{j=1}$.
Then (3.1) has a periodic solution if and only if
$$\d{\int^\pi_0}f(x)X_j(x)dx=0,\;\mbox{for}\;j=1,2,\cdots,m.$$Now, when $e(x)=0$, then $0$ is the first
eigenvalue of $L_0$ with $\sin x$ as a basis of the $0$-eigenspace.
Hence, (3.1) has a periodic solution if and only if
$$\d{\int^\pi_0}f(x)\sin xdx=0\;\mbox{or
}f(x)=\d{\sum^\infty_{j=2}}a_j\sin jx,\
\sum_{j=2}^{\infty}|a_j|^2<\infty.$$ Now suppose $e(x)\approx 0$.
The first eigenvalue $\lambda_e$ of $L_e$ is given by
$$\lambda_e=\d{\min_{\varphi\in
H^1_0(0,\pi),\|\varphi\|_{L^2(0,\pi)}=1}}J_e(\varphi,\varphi),$$where
$$J_e(\varphi,\varphi)=\d{\int^\pi_0}(\varphi^2_x-\varphi^2-e(x)\varphi^2)dx.$$
(See \cite{kn:[3]}). Hence,
$$\begin{array}{ll}
\lambda_e&\leq \d{\min_{\varphi\in
H^1_0(0,\pi),\|\varphi\|_{L^2(0,\pi)}=1}}\d{\int^\pi_0}(\varphi^2_x-\varphi^2)dx+\max|e(x)|\d{\int^\pi_0}\varphi^2dx\\
&\leq 0+\max|e(x)|\\
& \leq \max|e(x)|.
\end{array}\eqno{(3.2)}
$$
$$\lambda_e=J_e(\varphi_e,\varphi_e)=\d{\int^\pi_0}(\varphi_e)^2_xdx-\d{\int^\pi_0}(1+e(x))\varphi_e^2dx$$with
$\varphi_e$ the eigenvector corresponding to $\lambda_e$ and $\|\varphi_e\|_{L^2(0,\pi)}=1$.\medskip

Since $0$ is the first eigenvalue of $L_0$, we have
$$\begin{array}{ll}
\lambda_e&=\d{\int^\pi_0}((\varphi_e)^2_x-(\varphi_e)^2)dx-\d{\int^\pi_0}e(x)\varphi_e^2dx\\
&\geq -\max|e(x)|
\end{array}\eqno{(3.3)}
$$
By (3.2) and (3.3), we get $$|\lambda_e|\leq \max|e(x)|,\;\mbox{and
}\lambda_e\rightarrow 0 \;\mbox{as }e(x)\rightarrow 0.$$ Next,
consider the system with $e(x)+\lambda_e$ as the  perturbation in
the system conductivity:
$$\left\{\begin{array}{ll} y_t-y_{xx}-y-(e(x)+\lambda_e)y=f(x),\;\;\;\;\;\;&0\leq x\leq \pi,\;0\leq t\leq T,\\
y(0,t)=y(\pi,t)=0,&0\leq t\leq T.
\end{array}\right.\eqno{(3.4)}
$$
Then when $e(x)\approx 0$, we have $(e(x)+\lambda_e)\approx 0$. However, if (3.4) still has a periodic solution,
we have
$$
\d{\int^\pi_0}f(x)\varphi_edx=0.
$$
If this is the case for any given $f$, we then have
$$
\d{\int^\pi_0}\sin jx\varphi_e dx=0,\;\mbox{for }j=2,3,\cdots.
$$
This implies that $\varphi_e=C \sin x$ and thus
$$
-e(x)\sin x=\lambda_e \sin x,\;\mbox{or }e(x)=-\lambda_e.
$$
This is a contradiction unless $e(x)\equiv const.$. This shows that
for any non-constant small perturbation in $e(x)$, for most a priori
given $f$, the periodicity of the system will get lost.\bigskip

\section{Local stabilization}
\hspace*{0.5 cm}In this section, we consider the same problem as
studied in Section 2, but with the control only imposed over a
subregion $\omega\subset\Omega$ and time interval $E\subset [0,T]$
with $m(E)>0$.

For the proof of Theorem 1.2, we need the following lemma, whose
quantitative version in the Laplacian case can be found in
\cite{kn:[4]} and \cite{kn:[5]}:

\medskip
{\bf Lemma 4.1} {\it  Let $X_{ij}(\omega)=\d\int_\omega
X_i(x)X_j(x)dx$. Then the symmetric matrix
$X(\omega,k)=(X_{ij}(\omega))_{1\leq i,j\leq k}$ is positive
definite for any $k\geq 1$. In particular, it is
invertible.}\medskip

{\it Proof of Lemma 4.1:} Let $a=(a_1,a_2,\cdots,a_k)\in {\bf R^k}$
and let
$$
I(a,a)=\d\int_\omega|\sum^k_{j=1}a_jX_j(x)|^2dx.
$$
Then $$I(a,a)=a\cdot X(\omega,k)\cdot a^\tau,\mbox{ where
}a^\tau=\left(
\begin{array}{ccc}
 a_1
 \\
 a_2
\\
\vdots\\
 a_k
\end{array}
\right).
$$
Apparently, $I(a,a)\geq 0$. If $X(\omega,k)$ is not positive
definite, then there is a vector $a'=(a'_1,a'_2,\cdots,a'_k)\neq 0$
such that $I(a',a')=0$. Without loss of generality, assume that
$a'_k\not =0$. Hence,
$$\d\sum^k_{j=1}a'_jX_j(x)|_{\omega}=0.\eqno{(4.1)}$$ We thus get over $\omega$:
$$X_k(x)=\d\sum_{j<k}b_jX_j(x),\;\mbox{with }b_j=-\d\frac{a'_j}{a'_k}.\eqno{(4.2)}$$
Applying $(L_0)^m$ to (4.2) over $\omega$, we have
$$\lambda^m_k X_k(x)=\d\sum_{j<k}b_j\lambda^m_jX_j(x).$$
We get
$$X_k(x)=\d\sum_{j<k}b_j(\frac{\lambda_j}{\lambda_k})^mX_j(x)\mbox{ over }\omega.$$
Letting $m\rightarrow\infty$, we get over $\omega$
$$
X_k(x)=\d\sum_{k'\leq j<k}b_jX_j(x),\eqno{(4.3)}
$$
where
$$
 \left\{\begin{array}{ll}
\lambda_j=\lambda_k,\;\;\;&\mbox{for }j\geq k',\\
\lambda_j<\lambda_k,&\mbox{for }j<k'.
\end{array}\right.\eqno{(4.4)}
$$
By (4.4), we get over $\Omega$,
$$
L_0(X_k(x)-\d\sum_{k'\leq
j<k}b_jX_j(x))=\lambda_kX_k(x)-\d\sum_{k'\leq
j<l}b_j\lambda_jX_j(x)=\lambda_k[X_k(x)-\d\sum_{k'\leq
j<k}b_jX_j(x)].
$$
By (4.3) and the unique continuation for solutions of elliptic
equations, we get
$$
X_k(x)-\d\sum_{k'\leq j<k}b_jX_j(x)\equiv 0\;\mbox{over }\Omega.
$$
This contradicts the linear independence of the system
$\{X_j\}$.\endpf

\medskip
{\bf Proof of Theorem 1.2.}: Similar to the proof of Theorem 1.2, we
need only to find a vector $\vec{u}=(u_1,u_2,\cdots,u_{K_0})\in
{\mathbf R}^{K_0}$ such that
$$J^*_\omega(\vec{u})=-J_{0,\omega},$$
where
$$J^*_\omega(\vec{u})=(\langle v(x,T),X_1(x)\rangle,\langle v(x,T),X_2(x)\rangle,\cdots,\langle
v(x,T),X_{K_0}(x)\rangle)=(v_1(T),v_2(T),\cdots,v_{K_0}(T))
$$
with $v$ the solution of the following equation:
$$
\left\{\begin{array}{ll} v_t
+L_0v+e(x,t)v=\displaystyle{\sum_{j=1}^{K_0}}\chi_\omega(x)\chi_E(t)u_jX_j(x),
\;\;\;&
\mbox{in }\;\;Q,\\
 v=0, & \mbox{on }\;\; \Sigma,\\
v_j(0)=0, & \mbox{for }\;j\leq K_0,\\
v\in {\bf S_{K_0}}.
\end{array}\right.
\eqno{(4.5)}
$$
and
$$
J_{0,\omega}=((v_0)_1(T),(v_0)_2(T),\cdots,(v_0)_{K_0}(T))
$$
with $v_0$ the solution of the following system
$$
 \left\{\begin{array}{ll}
(v_0)_t +L_0v_0+e(x,t)v_0=-e(x,t)y_0, \;\;\;&
\mbox{in }\;\;Q,\\
 v_0=0, & \mbox{on }\;\; \Sigma,\\
(v_0)_{j}(0)=0, & \mbox{for }\;j\leq K_0,\\
v_0\in {\bf S_{K_0}}.
\end{array}\right.
\eqno{(4.6)}
$$
In the same way, if $J^*_\omega$ is not invertible, then for a
vector $\vec{\xi}=(\xi_1,\xi_2,\cdots,\xi_{K_0})$ with
$|\vec{\xi}|=1$, we have a solution to the following system:
$$
\left\{\begin{array}{ll} v_t
+L_0v+e(x,t)v=\displaystyle{\sum_{j=1}^{K_0}}\chi_\omega(x)\chi_E(t)\xi_jX_j(x),
\;\;\;&
\mbox{in }\;\;Q,\\
 v=0, & \mbox{on }\;\; \Sigma,\\
v_j(0)=0=v_j(T), & \mbox{for }\;j\leq K_0,\\
v\in {\bf S_{K_0}}.
\end{array}\right.
\eqno{}
$$
We then get
$$
v_j(t)'+\lambda_jv_j(t)+\d\int_\Omega
e(x,t)v(x,t)X_j(x)dx=\d\sum^{K_0}_{l=1}\xi_l\chi_E(t)X_{lj}(\omega),\;\mbox{for }j=1,2,\cdots,K_0.
$$
We similarly get
$$
(e^{\lambda_j t}v_j(t))'_t+e^{\lambda_j t}\d\int_\Omega e(x,t)v(x,t)X_j(x)dx=e^{\lambda_j
t}\d\sum^{K_0}_{l=1}\xi_l\chi_E(t)X_{lj}(\omega),\;\mbox{for }j=1,2,\cdots,K_0.
$$
$$
0+\d\int^T_0\int_\Omega e^{\lambda_j t}e(x,t)v(x,t)X_j(x)dxdt=\d\int^T_0e^{\lambda_j
t}\d\sum^{K_0}_{l=1}\xi_l\chi_E(t)X_{lj}(\omega)dt.
$$
We then get
$$
\left(
\begin{array}{ccc}
\d\int^T_0e^{\lambda_1 t}\chi_E(t)dt&\;&\;
 \\
 \;&\d\int^T_0e^{\lambda_2 t}\chi_E(t)dt&\;
\\
\;&\ddots&\;\\
\;&\;&\d\int^T_0e^{\lambda_{K_0} t}\chi_E(t)dt
\end{array}
\right)
X(\omega,K_0)\left(
\begin{array}{ccc}
 \xi_1
 \\
 \xi_2
\\
\vdots\\
 \xi_{K_0}
\end{array}
\right)$$
$$
=\left(
\begin{array}{ccc}
  \d\int^T_0\int_\Omega e^{\lambda_1 t}e(x,t)v(x,t)X_1(x)dxdt\\
\d\int^T_0\int_\Omega e^{\lambda_2 t}e(x,t)v(x,t)X_2(x)dxdt\\
\vdots\\
\d\int^T_0\int_\Omega e^{\lambda_{K_0} t}e(x,t)v(x,t)X_{K_0}(x)dxdt
\end{array}
\right)
$$
$$
\left(
\begin{array}{ccc}
 \xi_1
 \\
 \xi_2
\\
\vdots\\
 \xi_{K_0}
\end{array}
\right)=X(\omega,K_0)^{-1}
 \left(
\begin{array}{ccc}
(\d\int^T_0e^{\lambda_1 t}\chi_E(t)dt)^{-1}  \d\int_Q e^{\lambda_1 t}evX_1dxdt\\
(\d\int^T_0e^{\lambda_2 t}\chi_E(t)dt)^{-1} \d\int_Q e^{\lambda_2 t}evX_2dxdt\\
\vdots\\
(\d\int^T_0e^{\lambda_{K_0} t}\chi_E(t)dt)^{-1} \d\int_Q e^{\lambda_{K_0} t}evX_{K_0}dxdt
\end{array}
\right)\eqno{(4.7)}
$$
By Lemma 4.1, we know $X(\omega,K_0)^{-1}$ is a  bounded linear
operator from ${\mathbf R}^{K_0}$ to ${\mathbf R}^{K_0}$.

By the energy estimate in Theorem 2.1, we have for $v(x,t)$,
$$\begin{array}{ll}
\d\sup_{t\in [0,T]}\|v(\cdot,t)\|^2+\int^T_0\|\nabla
v(\cdot,t)\|^2dt&\leq C(system,K_0)\d\int_Q(\sum_{j=1}^{K_0}\chi_\omega(x)\chi_E(t)\xi_jX_j(x))^2dxdt\\
&\leq C(system,K_0)T|\vec{\xi}|^2\\
&\leq C(system,K_0).
\end{array}\eqno{(4.8)}$$
By the H$\ddot{o}$lder inequality, we have,
$$\begin{array}{ll}
\d\int_\Omega|evX_j|dx\leq
C(\Omega,N,q)\|e\|_{L^\infty(0,T;L^q(\Omega))}[1+\lambda_j^2+\|v(\cdot,t)\|^2+\|\nabla
v(\cdot,t)\|^2].
\end{array}\eqno{(4.9)}$$
Together with (4.8), we thus have
$$\begin{array}{ll}
\d\int^T_0\int_\Omega|evX_j|dxdt\leq
C(system,K_0)\|e\|_{L^\infty(0,T;L^q(\Omega))}.
\end{array}\eqno{(4.10)}$$
Back to (4.7), we have
$$\begin{array}{ll}  |\vec{\xi}|^2&\leq
C(system,\omega,K_0)\frac{1}{(m(E))^2}\|X(\omega,K_0)^{-1}\|^2\|e\|^2_{L^\infty(0,T;L^q(\Omega))}\\
&\leq
C(system,\omega,K_0)\frac{1}{(m(E))^2}\|e\|^2_{L^\infty(0,T;L^q(\Omega))}.
\end{array}$$
Hence, when $\|e\|^2_{L^\infty(0,T;L^q(\Omega))}$ is sufficient
small, we get $|\vec{\xi}|^2<1$. This gives a contradiction.
Therefore, we showed that $J^*_\omega$ is invertible under small
perturbation. By the same arguments as those in the proof of Theorem
1.1, we can also show the energy estimates as stated in Theorem 1.2.
 This completes the proof of Theorem 1.2. $\endpf$


\end{document}